\documentclass[11pt]{amsart}

\usepackage{amsmath,amssymb,amsthm,mathtools}
\usepackage[margin=1in]{geometry}
\usepackage{graphicx}
\usepackage{tikz}
\usepackage{enumitem}
\usepackage[colorlinks=true,citecolor=blue,linkcolor=blue,urlcolor=blue]{hyperref}

\theoremstyle{plain}
\newtheorem{theorem}{Theorem}[section]
\newtheorem{proposition}[theorem]{Proposition}
\newtheorem{lemma}[theorem]{Lemma}
\newtheorem{corollary}[theorem]{Corollary}

\theoremstyle{definition}
\newtheorem{definition}[theorem]{Definition}
\newtheorem{assumption}[theorem]{Assumption}
\newtheorem{example}[theorem]{Example}

\theoremstyle{remark}
\newtheorem{remark}[theorem]{Remark}

\newcommand{\Om}{\Omega}
\newcommand{\R}{\mathbb{R}}
\newcommand{\eps}{\varepsilon}
\newcommand{\Acal}{\mathcal A}

\DeclareMathOperator{\intr}{int}

\begin{document}

\title[ COINCIDENCE CRITERIA FOR PARTIALLY SEGREGATED SYSTEM]
{Coincidence criteria and junction obstructions for two partially segregated elliptic systems}

\author{Farid Bozorgnia}
\address{Department of Mathematics, New Uzbekistan University}
\email{f.bozorgnia@newuu.uz}

\date{August 2026}

\begin{abstract}
We compare two singularly perturbed elliptic systems with the same
partial-segregation constraint but different limiting mechanisms. System~A is
determined by harmonic-difference identities and admits an explicit
lower-envelope limit, whereas System~B converges, along subsequences, to a
constrained Dirichlet-energy minimizer. We show that the System~A profile need
not be stationary for the limiting variational problem. We establish
conditional coincidence criteria, give explicit data for which the interfaces
differ, and derive a planar triple-junction obstruction. In particular,
coincidence with the standard variational cone requires the triangle formed by
the harmonic gradients to be equilateral. This condition is necessary but not
sufficient.
\end{abstract}
\maketitle

\section{Introduction}\label{sec:intro}

Spatial segregation in elliptic systems arises in reaction--diffusion theory,
optimal partition problems, and harmonic maps into singular targets
\cite{CaffarelliLin,ContiTerraciniVerzini,WangZhang}. We consider three
nonnegative components in the partially segregated regime, where the limiting
constraint is
\[
u_1u_2u_3=0.
\]
Thus, pairwise coexistence is permitted, but all three components cannot be
positive at the same point.

We compare two singularly perturbed systems with the same boundary data and
the same limiting segregation constraint. The first, referred to as
System~A, is
\begin{equation}\label{eq:sysA}
\Delta u_i^\eps
=
\frac{1}{\eps}\prod_{j=1}^3 u_j^\eps
\quad\text{in }\Om,
\qquad
u_i^\eps=\phi_i
\quad\text{on }\partial\Om,
\qquad i=1,2,3.
\end{equation}
Since the reaction term in \eqref{eq:sysA} is independent of \(i\),
subtracting two equations shows that \(u_i^\eps-u_j^\eps\) is harmonic for
every \(\eps>0\). The common-reaction framework and its uniform regularity
theory were developed by Caffarelli and Roquejoffre
\cite{CaffarelliRoquejoffre}, while the explicit lower-envelope limit was
identified by Bozorgnia, Burger, and Fotouhi \cite{BBF}.  Recent extensions to systems with general \(k\)-wise interactions are
developed in~\cite{Giaretto}. A related planar \(\Gamma\)-convergence analysis of the penalized and
constrained formulations is given in
\cite{BozorgniaArakelyanGamma}. For comparison, refined triple-junction regularity in the fully segregated
harmonic-map and optimal-partition setting is developed in
\cite{OgnibeneVelichkovTriple}.

The second system, referred to as System~B, is
\begin{equation}\label{eq:sysB}
\Delta u_i^\eps
=
\frac{u_i^\eps}{\eps}
\prod_{j\ne i}(u_j^\eps)^2
\quad\text{in }\Om,
\qquad
u_i^\eps=\phi_i
\quad\text{on }\partial\Om,
\qquad i=1,2,3.
\end{equation}
It is the Euler--Lagrange system associated with
\begin{equation}\label{eq:Eeps-intro}
E_\eps(u)
=
\sum_{i=1}^3\int_\Om |\nabla u_i|^2\,dx
+
\frac{1}{\eps}\int_\Om (u_1u_2u_3)^2\,dx.
\end{equation}
Every sequence of minimizers with \(\eps\downarrow0\) admits a subsequence
converging to a minimizer of the Dirichlet energy under the partial-segregation
constraint; the corresponding variational theory was developed by Soave and
Terracini \cite{SoaveTerracini1,SoaveTerracini2}.

The two systems therefore select their limits through different mechanisms.
System~A is governed by global harmonic-difference identities and admits an
explicit lower-envelope representation. Its individual components may,
however, develop gradient jumps inside their positivity sets. By contrast,
System~B is selected through constrained Dirichlet-energy minimization and
satisfies the associated domain-stationarity identity. Consequently,
coincidence of the limiting profiles or interfaces is not automatic.

We establish sufficient conditions for coincidence, including regional
harmonic compatibility and boundary dominance with full support. We also give
explicit one-dimensional boundary data for which the limiting interfaces are
different. In dimension two, we show that the blow-up of the System~A
absent-component partition at a nondegenerate triple tie is the normal fan of
the triangle formed by the three harmonic gradients. If the same point is a
simple triple point of System~B and the interfaces coincide locally, this
gradient triangle must be equilateral. This condition is necessary but not
sufficient for local coincidence.

Section~\ref{sec:setting} establishes existence and compactness for
System~B and defines the selected limiting minimizer. Section~\ref{sec:limits}
derives the explicit System~A profile, records the structural properties of
selected System~B limits, and identifies a mechanism by which the System~A
profile fails to be stationary. Section~\ref{sec:coincidence} gives
coincidence and noncoincidence results. Section~\ref{sec:junctions} studies
junction geometry, and Section~\ref{sec:numerics} presents formal inner-layer
scales and numerical illustrations.

\section{Problem setting and compactness}\label{sec:setting}

Let \(\Om\subset\R^d\) be a bounded, smooth, connected domain, and fix
\(\alpha\in(0,1)\).

\begin{assumption}\label{ass:bdry}
The boundary data satisfy
\[
\phi_i\in C^{1,\alpha}(\partial\Om),
\qquad
\phi_i\ge0,
\qquad
\phi_1\phi_2\phi_3=0
\quad\text{on }\partial\Om,
\qquad i=1,2,3.
\]
\end{assumption}

The last condition excludes triple coexistence on the boundary but permits
two boundary data to be positive simultaneously.

For each \(i\), let \(h_i\) denote the harmonic extension of \(\phi_i\):
\[
\Delta h_i=0
\quad\text{in }\Om,
\qquad
h_i=\phi_i
\quad\text{on }\partial\Om.
\]
Set
\[
h_{ij}:=h_i-h_j.
\]
Then \(h_{ij}\) is harmonic with boundary trace \(\phi_i-\phi_j\), and
\[
h_{ij}=-h_{ji},
\qquad
h_{ik}=h_{ij}+h_{jk}.
\]

Let \(\gamma:H^1(\Om)\to H^{1/2}(\partial\Om)\) denote the trace operator.
We define
\begin{equation}\label{eq:classes}
\Acal
:=
\left\{
u=(u_1,u_2,u_3)\in H^1(\Om)^3:
u_i\ge0\ \text{a.e.},
\ \gamma u_i=\phi_i
\right\},
\end{equation}
and
\begin{equation}\label{eq:segregated-class}
S
:=
\left\{
u\in\Acal:
u_1u_2u_3=0\ \text{a.e. in }\Om
\right\}.
\end{equation}
The limiting Dirichlet energy is
\[
E_0(u):=\sum_{i=1}^3\int_\Om|\nabla u_i|^2\,dx.
\]

\begin{proposition}\label{prop:existence}
Let
\[
X
:=
\left\{
u\in H^1(\Om)^3:
\gamma u_i=\phi_i,\ i=1,2,3
\right\}.
\]
For every \(\eps>0\), the energy \(E_\eps\) admits a minimizer over \(X\).
Moreover,
\[
\min_XE_\eps=\min_{\Acal}E_\eps,
\]
and a minimizer may be chosen nonnegative. Every minimizer
\(u^\eps\in\Acal\) satisfies
\[
0\le u_i^\eps\le h_i
\qquad\text{a.e. in }\Om
\]
and is a weak solution of System~B.
\end{proposition}

\begin{proof}
Let \((u^{(n)})\subset X\) be a minimizing sequence. Choose
\[
M_i>\|\phi_i\|_{L^\infty(\partial\Om)}
\]
and define
\[
T_i(s):=\max\{-M_i,\min\{s,M_i\}\}.
\]
Componentwise truncation preserves the boundary trace and does not increase
either the Dirichlet term or the penalization term. Thus the minimizing
sequence may be assumed uniformly bounded in \(L^\infty(\Om)^3\).

After subtracting fixed \(H^1\)-extensions of the boundary data, Poincar\'e's
inequality and the Dirichlet term give a uniform \(H^1\)-bound. Passing to a
subsequence,
\[
u^{(n)}\rightharpoonup u
\quad\text{in }H^1(\Om)^3,
\qquad
u^{(n)}\to u
\quad\text{in }L^2(\Om)^3
\]
and almost everywhere. The uniform \(L^\infty\)-bound implies convergence of
the penalization term by dominated convergence, while the Dirichlet term is
weakly lower semicontinuous. Hence \(u\) minimizes \(E_\eps\) over \(X\).

For every \(v\in X\), the componentwise absolute value \(|v|\) belongs to
\(\Acal\) and satisfies
\[
E_\eps(|v|)=E_\eps(v).
\]
Since \(\Acal\subset X\), this proves
\[
\min_XE_\eps=\min_{\Acal}E_\eps
\]
and yields a nonnegative minimizer.

Let \(u^\eps\in\Acal\) be a minimizer. Replacing \(u_i^\eps\) by
\(\min\{u_i^\eps,h_i\}\) preserves the trace and nonnegativity, decreases the
penalization term, and does not increase the Dirichlet energy. Indeed, since
\(u_i^\eps-h_i\in H_0^1(\Om)\) and \(h_i\) is harmonic,
\[
\int_\Om|\nabla u_i^\eps|^2\,dx
=
\int_\Om|\nabla h_i|^2\,dx
+
\int_\Om|\nabla(u_i^\eps-h_i)|^2\,dx.
\]
 The truncation decreases the Dirichlet term by
\[
\int_\Om
\left|\nabla (u_i^\varepsilon-h_i)^+\right|^2\,dx
\]
and does not increase the penalization term. Minimality therefore yields
\((u_i^\varepsilon-h_i)^+=0\).

The resulting minimizer is bounded. Since it also minimizes over the affine
class \(X\), two-sided variations yield
\[
\int_\Om\nabla u_i^\eps\cdot\nabla\eta\,dx
+
\frac1\eps
\int_\Om
u_i^\eps\prod_{j\ne i}(u_j^\eps)^2\eta\,dx
=0
\]
for every \(\eta\in H_0^1(\Om)\) and \(i=1,2,3\).
\end{proof}

The following competitor proves that the constrained class is nonempty.

\begin{lemma}\label{lem:nonempty}
Under Assumption~\ref{ass:bdry}, the functions
\begin{equation}\label{eq:vcompetitor}
v_1:=\max\{h_{12},h_{13},0\},
\qquad
v_2:=v_1-h_{12},
\qquad
v_3:=v_1-h_{13}
\end{equation}
belong to \(S\).
\end{lemma}

\begin{proof}
The maximum of finitely many \(H^1\)-functions belongs to \(H^1\), and its
trace is the maximum of their traces. By construction,
\[
v_1\ge0,
\qquad
v_2=v_1-h_{12}\ge0,
\qquad
v_3=v_1-h_{13}\ge0.
\]
At every point, at least one of the three quantities \(0,h_{12},h_{13}\)
attains their maximum. Accordingly, at least one of \(v_1,v_2,v_3\)
vanishes, and therefore
\[
v_1v_2v_3=0
\quad\text{a.e. in }\Om.
\]

It remains to verify the trace. If \(\phi_1=0\) at a boundary point, then
\[
h_{12}=-\phi_2\le0,
\qquad
h_{13}=-\phi_3\le0,
\]
so \(v_1=0=\phi_1\). If \(\phi_1>0\), Assumption~\ref{ass:bdry} implies that
at least one of \(\phi_2,\phi_3\) vanishes. Hence one of \(h_{12},h_{13}\)
equals \(\phi_1\), while both are at most \(\phi_1\). Thus \(v_1=\phi_1\),
and the definitions of \(v_2,v_3\) give \(v_2=\phi_2\) and \(v_3=\phi_3\).
Therefore \(v\in S\).
\end{proof}

\begin{theorem}\label{thm:compactness}
Let \(\eps_n\downarrow0\), and let \(u^{\eps_n}\) minimize
\(E_{\eps_n}\) over \(\Acal\). After passing to a subsequence, there exists
\(u^B\in S\) such that
\[
u_i^{\eps_n}\to u_i^B
\qquad\text{strongly in }H^1(\Om),
\qquad i=1,2,3.
\]
The limit \(u^B\) minimizes \(E_0\) over \(S\), and
\[
E_{\eps_n}(u^{\eps_n})\to E_0(u^B),
\qquad
\frac1{\eps_n}
\int_\Om
(u_1^{\eps_n}u_2^{\eps_n}u_3^{\eps_n})^2\,dx
\to0.
\]
If the minimizer of \(E_0\) over \(S\) is unique, then the whole family
converges to it.
\end{theorem}

\begin{proof}
By Proposition~\ref{prop:existence},
\[
0\le u_i^\eps\le h_i.
\]
Let \(v\in S\) be the competitor from Lemma~\ref{lem:nonempty}. Since
\(v_1v_2v_3=0\),
\[
E_\eps(v)=E_0(v)=:C_0.
\]
Minimality gives
\[
E_\eps(u^\eps)\le C_0,
\]
and consequently
\[
\sum_{i=1}^3\int_\Om|\nabla u_i^\eps|^2\,dx\le C_0,
\qquad
\int_\Om(u_1^\eps u_2^\eps u_3^\eps)^2\,dx
\le C_0\eps.
\]
Together with the fixed traces and Poincar\'e's inequality, this yields a
uniform \(H^1\)-bound.

After passing to a subsequence,
\[
u^\eps\rightharpoonup u^B
\quad\text{in }H^1(\Om)^3,
\qquad
u^\eps\to u^B
\quad\text{in }L^2(\Om)^3.
\]
The trace and nonnegativity constraints pass to the limit. Moreover, the
uniform \(L^\infty\)-bounds imply
\[
\left\|
u_1^\eps u_2^\eps u_3^\eps-u_1^Bu_2^Bu_3^B
\right\|_{L^2}
\le
C\sum_{i=1}^3
\|u_i^\eps-u_i^B\|_{L^2}
\to0.
\]
The penalization estimate therefore gives
\[
u_1^Bu_2^Bu_3^B=0,
\]
so \(u^B\in S\).

For every \(w\in S\),
\[
E_\eps(u^\eps)\le E_\eps(w)=E_0(w).
\]
Hence
\[
E_0(u^B)
\le
\liminf_{\eps\to0}E_0(u^\eps)
\le
\limsup_{\eps\to0}E_\eps(u^\eps)
\le
E_0(w).
\]
Thus \(u^B\) minimizes \(E_0\) over \(S\). Taking \(w=u^B\) gives
\[
E_0(u^\eps)\to E_0(u^B)
\quad\text{and}\quad
E_\eps(u^\eps)\to E_0(u^B).
\]
The convergence of the Dirichlet norms, together with weak \(H^1\)
convergence, yields strong \(H^1\) convergence. Finally,
\[
\frac1\eps\int_\Om(u_1^\eps u_2^\eps u_3^\eps)^2\,dx
=
E_\eps(u^\eps)-E_0(u^\eps)
\to0.
\]
If the constrained minimizer is unique, every convergent subsequence has the
same limit, and therefore the entire family converges.
\end{proof}

The class \(S\) is nonconvex, and uniqueness of its Dirichlet-energy minimizer
is not assumed. Hence, throughout the remainder of the paper, \(u^B\) denotes
a fixed subsequential limit selected by a chosen sequence of penalized
minimizers.

\section{The two limits}\label{sec:limits}

\subsection{System A: the explicit lower-envelope profile}

 Subtracting two equations of System~A yields the following
harmonic-difference identity.
\begin{lemma}\label{lem:harmonic-diff-A}
For every $\eps>0$, every solution of System~A satisfies $u_i^\eps-u_j^\eps=h_{ij}$ in $\Om$. In the
limit, $u_i^A-u_j^A=h_{ij}$.
\end{lemma}

\begin{proof}
Subtracting the $j$-th equation from the $i$-th in \eqref{eq:sysA} gives $\Delta(u_i^\eps-u_j^\eps)=0$
with boundary trace $\phi_i-\phi_j$. Uniqueness of the harmonic extension gives the identity, which
passes to the limit.
\end{proof}

The following proposition turns the identity into the closed-form profile. Existence, convergence,
nonnegativity, and partial segregation of the System~A limit are imported from \cite{BBF};
Lemma~\ref{lem:harmonic-diff-A} then identifies the limit uniquely, so the full family converges.

\begin{proposition}\label{prop:A-profile}
The limit of System~A is
\begin{equation}\label{eq:Aprofile}
        u_i^A=h_i-\min_{1\le k\le3}h_k,
        \qquad i=1,2,3,
\end{equation}
which is nonnegative, partially segregated, locally Lipschitz in $\Om$, and satisfies
\[
        u_i^A-u_j^A=h_{ij},\qquad
        \{u_i^A=0\}=\{h_i=\textstyle\min_k h_k\},
        \qquad
        \{u_i^A>0\}=\{h_i>\textstyle\min_k h_k\}.
\]
\end{proposition}

\begin{proof}
By Lemma~\ref{lem:harmonic-diff-A}, $u_i^A-h_i$ is independent of $i$; write $u_i^A=h_i-c$.
Nonnegativity gives $c\le\min_k h_k$. Partial segregation gives $\min_i u_i^A=0$, so
$0=\min_i(h_i-c)=\min_i h_i-c$ and $c=\min_k h_k$. The support characterization is immediate, and
local Lipschitz continuity follows from smoothness of the $h_i$ and the Lipschitz character of the
minimum.
\end{proof}

Formula \eqref{eq:Aprofile} is equivalent to the representation in \cite{BBF}, since
$h_1-\min\{h_1,h_2,h_3\}=\max\{0,h_{12},h_{13}\}$ and then $u_2^A=u_1^A-h_{12}$, $u_3^A=u_1^A-h_{13}$.

For a partially segregated configuration
\(u=(u_1,u_2,u_3)\in S\), define
\[
Z_i(u):=\{x\in\Om:u_i(x)=0\},
\qquad
\omega_i(u):=\intr Z_i(u).
\]
We distinguish the absent-component interface
\[
\Gamma_{\mathrm{cell}}(u)
:=
\bigcup_{i=1}^3\bigl(\partial\omega_i(u)\cap\Om\bigr)
\]
from the positivity free boundary
\[
\Gamma_{\mathrm{pos}}(u)
:=
\bigcup_{i=1}^3
\bigl(\partial\{u_i>0\}\cap\Om\bigr).
\]
These sets may differ if a zero set contains a lower-dimensional component
with empty interior. For System~A, we write
\[
Z_i^A:=Z_i(u^A),
\qquad
\omega_i^A:=\omega_i(u^A).
\]

The following corollary describes the System~A partition through the absent-component cells and
records that the interfaces are selected relatively open portions of the nodal sets, not the full
nodal sets.
\begin{corollary}\label{cor:A-partition}
For distinct \(i,j,k\), set
\[
\Gamma_{ij}^A
:=
\partial\omega_i^A\cap\partial\omega_j^A\cap\Om.
\]
Then
\[
\Gamma_{ij}^A
\subset
\{h_i=h_j\le h_k\}.
\]
At every point where exactly \(h_i\) and \(h_j\) attain the minimum, the
interface is locally contained in the relatively open nodal portion
\[
\{h_i=h_j<h_k\}\subset\{h_{ij}=0\}.
\]
Away from all tie sets, the cells are characterized by
\begin{center}
\begin{tabular}{c|c}
Cell & Sign conditions\\ \hline
\(\omega_3^A\) & \(h_{13}>0,\ h_{23}>0\)\\
\(\omega_2^A\) & \(h_{12}>0,\ h_{23}<0\)\\
\(\omega_1^A\) & \(h_{12}<0,\ h_{13}<0\)
\end{tabular}
\end{center}
\end{corollary}

\begin{proof}
On $\partial\omega_i^A\cap\partial\omega_j^A$ both $u_i^A$ and $u_j^A$ vanish by continuity, so
$h_{ij}=u_i^A-u_j^A=0$ and $h_i=h_j$; this value is the minimum, so $h_i=h_j\le h_k$. The sign table
follows from \eqref{eq:Aprofile} and $h_{13}=h_{12}+h_{23}$.
\end{proof}

\begin{proposition}\label{prop:transversality}
Assume $0$ is a regular value of each $h_{ij}$ in $\Om$, that is $\nabla h_{ij}\ne0$ on
$\{h_{ij}=0\}\cap\Om$. Then each nodal set $\{h_{ij}=0\}$ is a smooth hypersurface, and the selected
System~A interfaces are smooth away from triple ties.
\end{proposition}

\begin{proof}
This is the implicit function theorem at regular values. The selected interface is the relatively
open portion $\{h_i=h_j<h_k\}$ of $\{h_{ij}=0\}$, whose relative closure may contain triple ties, and
which is relatively open in the smooth hypersurface.
\end{proof}

 The next proposition identifies a nondegenerate internal-exchange mechanism
under which the System~A profile fails to be stationary for the constrained
Dirichlet energy.

\begin{proposition}\label{prop:A-not-stationary}
Let $x_0$ be an interior point of $\{u_i^A>0\}$ at which the minimum is exchanged between the other
two extensions, $h_j(x_0)=h_k(x_0)<h_i(x_0)$, with $\nabla(h_j-h_k)(x_0)\ne0$. Set $f:=h_j-h_k$.
There is a neighborhood $U$ of $x_0$ on which $u_i^A>0$ and
\begin{equation}\label{eq:singular-laplacian}
        \Delta u_i^A=|\nabla(h_j-h_k)|\,\mathcal H^{d-1}\!\restriction\big(\{h_j=h_k\}\cap U\big).
\end{equation}
In particular, $\Delta u_i^A$ has a nonzero singular part inside $\{u_i^A>0\}$, and $u^A$ is not a
critical point of $E_0$ over $S$.
\end{proposition}

\begin{proof}
Since $h_j(x_0)=h_k(x_0)<h_i(x_0)$ and $\nabla f(x_0)\ne0$, choose a neighborhood $U$ of $x_0$ on
which $h_i>\max\{h_j,h_k\}$, so $u_i^A>0$, and on which $\{f=0\}$ is a smooth hypersurface with
$\nabla f\ne0$. On $U$ the minimum is attained by $h_j$ or $h_k$, so $u_i^A=h_i-\min\{h_j,h_k\}$ and
$\min\{h_j,h_k\}=\tfrac12(h_j+h_k-|f|)$. The function $|f|$ has a jump of $2|\nabla f|$ in its normal
derivative across $\{f=0\}$, so
\[
\Delta|f|=2|\nabla f|\,\mathcal H^{d-1}\!\restriction(\{f=0\}\cap U).
\]
Therefore
\[
        \Delta u_i^A=\Delta h_i-\tfrac12(\Delta h_j+\Delta h_k)+\tfrac12\Delta|f|
        =|\nabla f|\,\mathcal H^{d-1}\!\restriction\big(\{f=0\}\cap U\big),
\]
which is \eqref{eq:singular-laplacian}. Since $u_i^A>0$ on $U$ and at least one other component
vanishes there, compactly supported variations of $u_i$ alone are admissible, and a stationary point
would satisfy $\Delta u_i^A=0$ on $U$. This fails.
\end{proof}

\subsection{System B: the variational limit}

 The following properties follow from 
\cite{SoaveTerracini1,SoaveTerracini2}. The domain-stationarity identity is
variational information absent from System~A.
\begin{theorem}[Structural properties of System~B limits]
\label{thm:ST}
Let \(u^B=(u_1^B,u_2^B,u_3^B)\in S\) be a nontrivial selected limit of
minimizers of \(E_\eps\). Then:

\begin{enumerate}[label=\textup{(\roman*)}]

\item
\[
u^B\in C_{\mathrm{loc}}^{0,3/4}(\Om;\R^3),
\]
and the exponent \(3/4\) is optimal for this class.

\item
Each component is harmonic in its positivity set,
\[
\Delta u_i^B=0
\qquad\text{in }\{u_i^B>0\},
\]
and
\[
u_1^Bu_2^Bu_3^B=0
\qquad\text{in }\Om.
\]

\item
The configuration \(u^B\) is stationary under compactly supported domain
variations: for every \(Y\in C_c^1(\Om;\R^d)\),
\begin{equation}
\label{eq:domain-stationarity}
\int_\Om
\left[
\left(\sum_{i=1}^3|\nabla u_i^B|^2\right)\operatorname{div} Y
-
2\sum_{i=1}^3
(DY\,\nabla u_i^B)\cdot\nabla u_i^B
\right]dx
=0.
\end{equation}
Consequently, for every \(x_0\in\Om\) and almost every \(r>0\) with
\(\overline{B_r(x_0)}\subset\Om\),
\begin{equation}
\label{eq:local-pohozaev}
\int_{\partial B_r(x_0)}
\sum_{i=1}^3|\nabla u_i^B|^2\,d\sigma
=
\frac{d-2}{r}
\int_{B_r(x_0)}
\sum_{i=1}^3|\nabla u_i^B|^2\,dx
+
2\int_{\partial B_r(x_0)}
\sum_{i=1}^3(\partial_\nu u_i^B)^2\,d\sigma .
\end{equation}

\item
Either two components vanish identically and the remaining component is
harmonic in \(\Om\), or the positivity free boundary
\(\Gamma_{\mathrm{pos}}(u^B)\) is locally a finite union of smooth
hypersurfaces outside a singular set. In dimension \(d=2\), the singular
set is discrete. In the latter case, with
\[
Z_i^B:=\{u_i^B=0\},
\qquad
\omega_i^B:=\intr Z_i^B,
\]
one has
\[
Z_i^B=\overline{\omega_i^B},
\qquad
\Om=\bigcup_{i=1}^3\overline{\omega_i^B}.
\]

\end{enumerate}
\end{theorem}

 The optimal regularity and its sharpness in \textup{(i)} follow from
\cite[Theorems~1.1, 1.6, and 1.7]{SoaveTerracini2}.
Harmonicity, partial segregation, and the local Pohozaev identity follow
from \cite[Theorem~1.7]{SoaveTerracini1} and
\cite[(1.9) and (2.5)]{SoaveTerracini2}.
The stress-energy identity in \textup{(iii)} also follows directly from
stationarity of the constrained minimizer under compactly supported domain
variations. The free-boundary and cell structure in \textup{(iv)} follow
from \cite[Theorems~1.1 and 1.4 and Remark~1.5]{SoaveTerracini2}.
Existence with fixed trace is given by
\cite[Theorem~1.2]{SoaveTerracini1}.

Local minimality
implies the harmonic-difference property stated in the next lemma.
\begin{lemma}\label{lem:diff-harmonic-B}
Let $\{i,j,k\}=\{1,2,3\}$ and $G_k:=\{u_k^B>0\}$. Then $u_i^Bu_j^B=0$ in $G_k$, and
$w_{ij}:=u_i^B-u_j^B$ is harmonic in $G_k$.
\end{lemma}

\begin{proof}
Since $u_k^B>0$ in $G_k$ and $u_1^Bu_2^Bu_3^B=0$, we have $u_i^Bu_j^B=0$ in $G_k$. Fix $D\Subset G_k$
and $\eta\in C_c^\infty(D)$. For real $t$ set $\widetilde w=w_{ij}+t\eta$, $\widetilde u_i=\widetilde w^+$,
$\widetilde u_j=\widetilde w^-$, and $\widetilde u_k=u_k^B$ in $D$, unchanged outside $D$. This is
admissible, since $\widetilde u_i$ and $\widetilde u_j$ have disjoint supports and the trace on
$\partial D$ is preserved. Because positive and negative parts have disjoint supports,
$|\nabla\widetilde u_i|^2+|\nabla\widetilde u_j|^2=|\nabla\widetilde w|^2$ a.e., and similarly for
$w_{ij}$. Local minimality gives $\int_D|\nabla w_{ij}|^2\le\int_D|\nabla(w_{ij}+t\eta)|^2$ for all
$t$, so differentiating at $t=0$ yields $\int_D\nabla w_{ij}\cdot\nabla\eta=0$. Hence $w_{ij}$ is
harmonic in $G_k$.
\end{proof}

We use the following blow-up-based notion of a simple triple point.
\begin{definition}\label{def:simple}
Let $d=2$. A point $p\in\Om$ is a simple triple point of $u^B$ if some blow-up of $u^B$ at $p$
equals, up to rotation, scaling, and permutation of components, the standard $3/4$-homogeneous
partially segregated minimizing cone of \cite[Theorem 1.7]{SoaveTerracini2}.
\end{definition}

\section{Coincidence criteria}\label{sec:coincidence}

System~A is determined by global harmonic-difference identities, whereas
System~B is selected by constrained minimization of \(E_0\).
Proposition~\ref{prop:A-not-stationary} shows that the System~A profile need
not be stationary for \(E_0\); coincidence therefore requires additional
compatibility.

\begin{example}[Coincident interface, distinct profiles]\label{ex:oneD}
Let \(\Om=(0,1)\) and prescribe
\[
(\phi_1(0),\phi_2(0),\phi_3(0))=(1,0,1),
\qquad
(\phi_1(1),\phi_2(1),\phi_3(1))=(0,1,1).
\]
Then
\[
h_1=1-x,\qquad h_2=x,\qquad h_3=1,
\]
and
\[
u_1^A=(1-2x)^+,\qquad
u_2^A=(2x-1)^+,\qquad
u_3^A=1-\min\{x,1-x\}.
\]
The unique System~B minimizer is
\[
u^B=\bigl((1-2x)^+,(2x-1)^+,1\bigr).
\]
Thus both systems have interface \(\{1/2\}\), but
\[
E_0(u^A)=5,
\qquad
E_0(u^B)=4,
\]
so their limiting profiles differ.
\end{example}

\begin{proof}
Let \(u\in S\) have the prescribed traces. If \(u_3>0\) on \((0,1)\),
then \(u_1u_2=0\). Setting \(w=u_1-u_2\), we have
\(u_1=w^+\), \(u_2=w^-\), \(w(0)=1\), and \(w(1)=-1\). Hence
\[
\int_0^1\bigl(|u_1'|^2+|u_2'|^2\bigr)\,dx
=
\int_0^1|w'|^2\,dx
\ge4,
\]
with equality for \(w=1-2x\); moreover, \(u_3\equiv1\) minimizes its
Dirichlet energy.

If \(u_3(x_0)=0\) for some \(x_0\in(0,1)\), then
\[
\int_0^1|u_3'|^2\,dx
\ge
\frac1{x_0}+\frac1{1-x_0}
\ge4,
\]
while the endpoint data imply
\[
\int_0^1|u_1'|^2\,dx\ge1,
\qquad
\int_0^1|u_2'|^2\,dx\ge1.
\]
Such a competitor has energy at least \(6\), proving the claim.
\end{proof}

The following regional compatibility condition identifies the
absent-component cells.

\begin{proposition}\label{prop:regional}
Let \(u^B\) be a selected minimizer, assume
\[
h_{ij}\not\equiv0
\qquad(i\ne j),
\]
and suppose that, for every set of distinct indices
\(\{i,j,k\}=\{1,2,3\}\),
\begin{equation}\label{eq:regional-hyp}
u_i^B-u_j^B=h_{ij}
\qquad\text{in }G_k:=\{u_k^B>0\}.
\end{equation}
Then
\[
\omega_i^B=\omega_i^A
\qquad(i=1,2,3),
\]
and consequently
\[
\Gamma_{\mathrm{cell}}(u^B)
=
\Gamma_{\mathrm{cell}}(u^A).
\]
\end{proposition}

\begin{proof}
The nondegeneracy assumption excludes the exceptional alternative in
Theorem~\ref{thm:ST}\textup{(iv)} and ensures
\[
\omega_j^A=\{h_j<\min\{h_i,h_k\}\}.
\]

Let \(x\in\omega_j^B\). Since the open cells are pairwise disjoint and
\(Z_\ell^B=\overline{\omega_\ell^B}\), one has
\(x\in G_i\cap G_k\). Applying \eqref{eq:regional-hyp} in \(G_k\) and
\(G_i\) gives
\[
h_i(x)>h_j(x),
\qquad
h_k(x)>h_j(x),
\]
so \(x\in\omega_j^A\). Hence \(\omega_j^B\subset\omega_j^A\).

Conversely, if \(x\in\omega_j^A\), choose
\(B_r(x)\Subset\omega_j^A\). The first inclusion gives
\[
\omega_i^B\subset\omega_i^A,
\qquad
\omega_k^B\subset\omega_k^A,
\]
so \(B_r(x)\) is disjoint from
\(\overline{\omega_i^B}\cup\overline{\omega_k^B}\). Since
\[
\Om=\bigcup_{\ell=1}^3\overline{\omega_\ell^B},
\]
we obtain
\[
B_r(x)\subset Z_j^B,
\]
and therefore \(x\in\omega_j^B\). Taking relative boundaries proves the
interface identity.
\end{proof}

\begin{corollary}\label{cor:global}
If
\[
u_i^B-u_j^B=h_{ij}
\qquad\text{in }\Om
\]
for every \(i,j\), then \(u^B=u^A\). Consequently, both their
absent-component interfaces and positivity free boundaries coincide.
\end{corollary}

\begin{proof}
The functions \(u_i^B-h_i\) are independent of \(i\), so
\(u_i^B=h_i-c\) for a common function \(c\). Nonnegativity and partial
segregation give
\[
0=\min_i u_i^B=\min_i h_i-c,
\]
hence \(c=\min_i h_i\) and \(u^B=u^A\).
\end{proof}

The next result gives a boundary-dominance criterion, conditional on full
support of the dominant System~B component.

\begin{theorem}\label{thm:dominated}
Let \(\{i,j,k\}=\{1,2,3\}\) and assume
\begin{equation}\label{eq:dominance}
\phi_k\ge\max\{\phi_i,\phi_j\}
\quad\text{on }\partial\Om,
\qquad
\phi_k\not\equiv\phi_i,
\qquad
\phi_k\not\equiv\phi_j.
\end{equation}
Then
\[
h_k>\max\{h_i,h_j\}
\quad\text{in }\Om,
\qquad
\phi_i\phi_j=0
\quad\text{on }\partial\Om,
\]
and
\[
u_i^A=h_{ij}^+,\qquad
u_j^A=h_{ij}^-,\qquad
u_k^A=h_k-\min\{h_i,h_j\}>0.
\]
If \(u_k^B>0\) in \(\Om\), then
\[
u_k^B=h_k,\qquad
u_i^B=h_{ij}^+,\qquad
u_j^B=h_{ij}^-,
\]
and
\[
\Gamma_{\mathrm{pos}}(u^A)
=
\Gamma_{\mathrm{pos}}(u^B).
\]
Moreover,
\[
u_k^B-u_k^A=\min\{h_i,h_j\}.
\]
In particular, the dominant profiles agree exactly when
\(\min\{h_i,h_j\}\equiv0\), and if both nondominant boundary data are
nontrivial, then \(u_k^B>u_k^A\) in \(\Om\).
\end{theorem}

\begin{proof}
The harmonic functions \(h_k-h_i\) and \(h_k-h_j\) have nonnegative,
nontrivial boundary traces. The strong maximum principle gives
\[
h_k>h_i,\qquad h_k>h_j
\quad\text{in }\Om.
\]
If \(\phi_i\) and \(\phi_j\) were simultaneously positive, then
\eqref{eq:dominance} would imply \(\phi_k>0\), contradicting the boundary
segregation condition. Thus \(\phi_i\phi_j=0\).

Since
\[
\min_m h_m=\min\{h_i,h_j\},
\]
formula \eqref{eq:Aprofile} gives the System~A profile. If \(u_k^B>0\) in
\(\Om\), then \(u_k^B\) is harmonic with trace \(\phi_k\), so \(u_k^B=h_k\).
Lemma~\ref{lem:diff-harmonic-B} similarly gives
\[
u_i^B-u_j^B=h_{ij}.
\]
Together with \(u_i^Bu_j^B=0\), this yields
\[
u_i^B=h_{ij}^+,
\qquad
u_j^B=h_{ij}^-.
\]
The stated interface and profile identities follow immediately.
\end{proof}

\subsection{A one-dimensional noncoincidence example}
\label{sec:oned}

\begin{example}\label{ex:1d-noncoincidence}
Let \(\Om=(0,1)\) with boundary values
\[
(\phi_1,\phi_2,\phi_3)(0)=(3,0,1),
\qquad
(\phi_1,\phi_2,\phi_3)(1)=(0,3,1).
\]
Then
\[
h_1=3(1-x),\qquad h_2=3x,\qquad h_3=1.
\]
System~A has interface \(\{1/3,2/3\}\) and energy \(E_0(u^A)=36\).
The unique System~B minimizer is
\[
u_1^B=(3-4x)^+,\qquad
u_2^B=(4x-1)^+,\qquad
u_3^B=(1-4x)^++(4x-3)^+,
\]
with interface \(\{1/4,3/4\}\) and energy \(E_0(u^B)=32\).

\begin{proof}
For any \(u\in S\), set
\[
x_1:=\min\{u_1=0\},
\qquad
x_2:=\max\{u_2=0\}.
\]
Then
\[
\int_0^1|u_1'|^2\,dx\ge\frac9{x_1},
\qquad
\int_0^1|u_2'|^2\,dx\ge\frac9{1-x_2}.
\]
If \(x_1\le x_2\), then \(E_0(u)\ge36\). If \(x_2<x_1\), partial
segregation gives \(u_3=0\) on \((x_2,x_1)\), and hence
\[
E_0(u)\ge
\frac9{x_1}+\frac1{1-x_1}
+\frac1{x_2}+\frac9{1-x_2}
\ge32.
\]
Equality holds only for \(x_1=3/4\), \(x_2=1/4\), and determines the
displayed piecewise-affine profile.
\end{proof}
\end{example}

Thus the two interfaces may differ even in one dimension.
\section{Junction geometry}\label{sec:junctions}

We study the blow-up of the System~A absent-component partition near a
planar triple tie \(p\), where
\[
h_1(p)=h_2(p)=h_3(p).
\]
Set
\[
a_i:=\nabla h_i(p),
\qquad
\mathcal T_p:=\operatorname{conv}\{a_1,a_2,a_3\}.
\]
When \(a_1,a_2,a_3\) are not collinear, their first-order expansions
determine the tangent partition. We use the minimizing-normal convention
\[
C_i:=\{x\in\R^2:a_i\cdot x\le a_j\cdot x
\text{ for all }j\},
\]
so that \(\{C_1,C_2,C_3\}\) is the minimizing normal fan of
\(\mathcal T_p\).

\begin{theorem}\label{thm:gradient-triangle}
Let \(d=2\), and let \(p\in\Om\) be a triple tie with
\(a_1,a_2,a_3\) not collinear. Then the blow-up of the System~A
absent-component partition at \(p\) is the minimizing normal fan of
\(\mathcal T_p\). More precisely, the tangent cone of
\(\overline{\omega_i^A}\) is
\[
C_i
=
\{x\in\R^2:(a_i-a_j)\cdot x\le0
\text{ for all }j\}.
\]
If \(\alpha_i\) is the interior angle of \(\mathcal T_p\) at \(a_i\),
then \(C_i\) has opening \(\pi-\alpha_i\). The three interface rays are
orthogonal to the corresponding gradients \(\nabla h_{ij}(p)\), and they
are separated by angles \(2\pi/3\) if and only if
\[
|\nabla h_{12}(p)|
=
|\nabla h_{23}(p)|
=
|\nabla h_{31}(p)|.
\]
\end{theorem}

\begin{proof}
For \(x\) in a fixed compact set,
\[
h_i(p+rx)
=
h_i(p)+r\,a_i\cdot x+O(r^2)
\qquad(r\downarrow0).
\]
Hence, away from the boundaries of the cones,
\[
h_i(p+rx)<h_j(p+rx)\ \text{for all }j\ne i
\]
eventually holds exactly when
\[
a_i\cdot x<a_j\cdot x\ \text{for all }j\ne i.
\]
Thus the rescaled cells converge locally to the cones \(C_i\).

For \(\{i,j,k\}=\{1,2,3\}\),
\[
C_i
=
\{x:(a_j-a_i)\cdot x\ge0,\ 
      (a_k-a_i)\cdot x\ge0\}.
\]
Its boundary rays are therefore orthogonal to
\(a_j-a_i=\nabla h_{ji}(p)\) and
\(a_k-a_i=\nabla h_{ki}(p)\). Since the angle between these two vectors
is \(\alpha_i\), the opening of \(C_i\) is \(\pi-\alpha_i\).
All three openings equal \(2\pi/3\) precisely when
\(\alpha_1=\alpha_2=\alpha_3=\pi/3\), that is, when
\(\mathcal T_p\) is equilateral. Its side lengths are
\[
|a_i-a_j|=|\nabla h_{ij}(p)|,
\]
which proves the final equivalence.
\end{proof}
\begin{corollary}\label{cor:equilateral}
Let \(d=2\), and let \(p\in\Om\) satisfy
\[
h_1(p)=h_2(p)=h_3(p),
\]
with
\[
\nabla h_1(p),\ \nabla h_2(p),\ \nabla h_3(p)
\]
not collinear. Assume that \(p\) is a simple triple point of System~B and
that, for some neighborhood \(U\) of \(p\),
\[
\Gamma_{\mathrm{cell}}(u^A)\cap U
=
\Gamma_{\mathrm{pos}}(u^B)\cap U.
\]
Then
\[
|\nabla h_{12}(p)|
=
|\nabla h_{23}(p)|
=
|\nabla h_{31}(p)|.
\]
\end{corollary}

\begin{proof}
Theorem~\ref{thm:gradient-triangle} gives a unique tangent fan for the
System~A partition. Local coincidence of the interfaces implies that this
fan is also the interface set of the standard System~B blow-up at \(p\).
The latter is an equiangular tripod, so
Theorem~\ref{thm:gradient-triangle} yields the stated equality.
\end{proof}

\begin{figure}[t]
\centering
\begin{tikzpicture}[scale=1.0]
\begin{scope}[xshift=-3.4cm]
  \node at (0,2.3) {System A};
  \draw[thick] (0,0)--(1.6,0.55) node[right] {$r_{23}$};
  \draw[thick] (0,0)--(-1.5,0.9) node[left] {$r_{12}$};
  \draw[thick] (0,0)--(-0.2,-1.6) node[below] {$r_{13}$};
  \fill (0,0) circle (1.5pt) node[above right] {$p$};
  \node at (-0.55,0.55) {$\omega_3^A$};
  \node at (0.55,-0.35) {$\omega_2^A$};
  \node at (-1.0,-0.55) {$\omega_1^A$};
\end{scope}
\begin{scope}[xshift=3.4cm]
  \node at (0,2.3) {System B};
  \draw[thick] (0,0)--(0,1.65);
  \draw[thick] (0,0)--({1.65*cos(-30)},{1.65*sin(-30)});
  \draw[thick] (0,0)--({1.65*cos(210)},{1.65*sin(210)});
  \fill (0,0) circle (1.5pt) node[above right] {$p$};
  \node at (0.5,0.55) {$120^\circ$};
  \node at (0.55,-0.6) {$120^\circ$};
  \node at (-0.75,-0.15) {$120^\circ$};
\end{scope}
\end{tikzpicture}
\caption{Planar triple-junction obstruction. The System~A interface rays
form the minimizing normal fan of the gradient triangle \(\mathcal T_p\),
whereas the standard System~B blow-up is an equiangular tripod. Local
coincidence therefore requires \(\mathcal T_p\) to be equilateral.}
\label{fig:triple-obstruction}
\end{figure}

\begin{remark}\label{rem:limitations}
The equilateral condition is only necessary. It determines the opening
angles of the tangent sectors, but not their labeled orientation or the
higher-order geometry of the interfaces. Additional symmetry or uniqueness
assumptions would be required to align the labeled rays, and even equality
of tangent cones does not imply equality of the interface curves.
\end{remark}
\subsection{Higher-dimensional extension}

\begin{proposition}\label{prop:higherdim}
Let \(d\ge3\), and let \(p\) lie on a smooth codimension-two common
junction stratum \(\Sigma\), with normal plane
\[
N_p:=(T_p\Sigma)^\perp.
\]
Assume that the System~B blow-up at \(p\) is the cylindrical extension of
the standard planar cone and that
\(\Gamma_{\mathrm{cell}}(u^A)\) and
\(\Gamma_{\mathrm{pos}}(u^B)\) coincide near \(p\). Then
\[
P_{N_p}\nabla h_1(p),\qquad
P_{N_p}\nabla h_2(p),\qquad
P_{N_p}\nabla h_3(p)
\]
form an equilateral triangle in \(N_p\).
\end{proposition}

\begin{proof}
Restrict the first-order expansions of the \(h_i\) to the transverse plane
\(p+N_p\). The System~A tangent partition is the minimizing normal fan of
the projected-gradient triangle, while the transverse System~B blow-up is
an equiangular tripod. Coincidence gives the conclusion.
\end{proof}

\section{Numerical illustrations}\label{sec:numerics}

Let
\[
\Om=(-1,1)^2,
\qquad
h=\frac{2}{N-1},
\]
and discretize the Laplacian on a uniform \(N\times N\) grid by the standard
five-point stencil. Writing
\(\theta=\operatorname{atan2}(y,x)\), we prescribe
\begin{equation}\label{eq:caps}
\phi_i(\theta)
=
\max\{0,\cos(\theta-\theta_i)\},
\qquad i=1,2,3,
\end{equation}
where angular differences are understood modulo \(2\pi\). For both choices
of \((\theta_1,\theta_2,\theta_3)\) below, the three positive caps have empty
common intersection, and hence
\(\phi_1\phi_2\phi_3=0\) on \(\partial\Om\).

System~A is obtained by solving the three discrete harmonic-extension
problems and applying \eqref{eq:Aprofile}. For System~B, we use continuation
in
\[
\eps\in\{10^{-2},10^{-3},10^{-4},3\cdot10^{-5},10^{-5}\}
\]
and a block Gauss--Seidel iteration. With the other components frozen, each
update solves
\begin{equation}\label{eq:block}
\left(
-\Delta_h
+
\frac{1}{\eps}
(u_j^{\mathrm{latest}})^2
(u_k^{\mathrm{latest}})^2
\right)
u_i^{\mathrm{new}}
=0
\end{equation}
with the prescribed Dirichlet data. The first continuation stage is
initialized by the System~A profile. Since the penalized energy is
nonconvex, the resulting iterate is only a computed critical-point candidate
and is not guaranteed to be the global minimizer.

At finite \(\eps\), all three components are positive in the interior.
We therefore compare proxy partitions. At each grid node \(x\), define
\[
\ell_A^h(x)
:=
\min\operatorname*{arg\,min}_{1\le i\le3}h_i^h(x),
\qquad
\ell_B^{\eps,h}(x)
:=
\min\operatorname*{arg\,min}_{1\le i\le3}u_i^{\eps,h}(x),
\]
where the smallest index breaks ties. The corresponding proxy interfaces
are represented by the midpoints of grid edges across which the label
changes. Their symmetric Hausdorff distance is computed in
\[
Q_{0.85}:=[-0.85,0.85]^2
\]
to reduce boundary effects.

The System~A triple tie is approximated by the interior grid point
\[
p_h
\in
\operatorname*{arg\,min}_{x}
\sum_{i<j}\bigl(h_i^h(x)-h_j^h(x)\bigr)^2.
\]
Central differences at \(p_h\) provide the discrete gradient triangle
\[
\mathcal T_{p_h}^h
:=
\operatorname{conv}
\{\nabla_hh_1(p_h),\nabla_hh_2(p_h),\nabla_hh_3(p_h)\}.
\]

A formal balance suggests an inner scale \(\eps^{1/4}\) near a regular
double interface, assuming linear vanishing of the switching components,
and a scale \(\eps^{1/5}\) near a standard \(3/4\)-homogeneous triple core.
For \(N=121\) and \(\eps=10^{-5}\), these scales correspond to approximately
\(3.4h\) and \(6h\), respectively, so the layers are only marginally
resolved.

We consider
\[
\begin{aligned}
\text{Datum I:}\quad&
(\theta_1,\theta_2,\theta_3)
=
(0,2\pi/3,4\pi/3),\\
\text{Datum II:}\quad&
(\theta_1,\theta_2,\theta_3)
=
(0,1.3,3.3).
\end{aligned}
\]
 Although the boundary peaks are separated by \(2\pi/3\), the square is not
invariant under this rotation; hence Datum~I has no exact threefold symmetry. It is,
however, invariant under reflection across the \(x\)-axis combined with the
exchange of components \(2\) and \(3\). The computed tie lies on this axis,
and the corresponding two side lengths agree to the displayed precision.

For Datum~I, the discrete gradient triangle is nearly equilateral and the
two proxy interfaces differ by approximately one mesh width. This is
consistent with, but does not prove, coincidence of the limiting interfaces.
For Datum~II, the gradient triangle is far from equilateral, providing strong
numerical evidence that the necessary condition in
Corollary~\ref{cor:equilateral} fails. The observed displacement of the
finite-\(\eps\) proxy is consistent with this obstruction, but does not prove
the existence or location of a System~B triple point.
\begin{figure}[h]
\centering
\IfFileExists{I6.png}{\includegraphics[width=.9\textwidth]{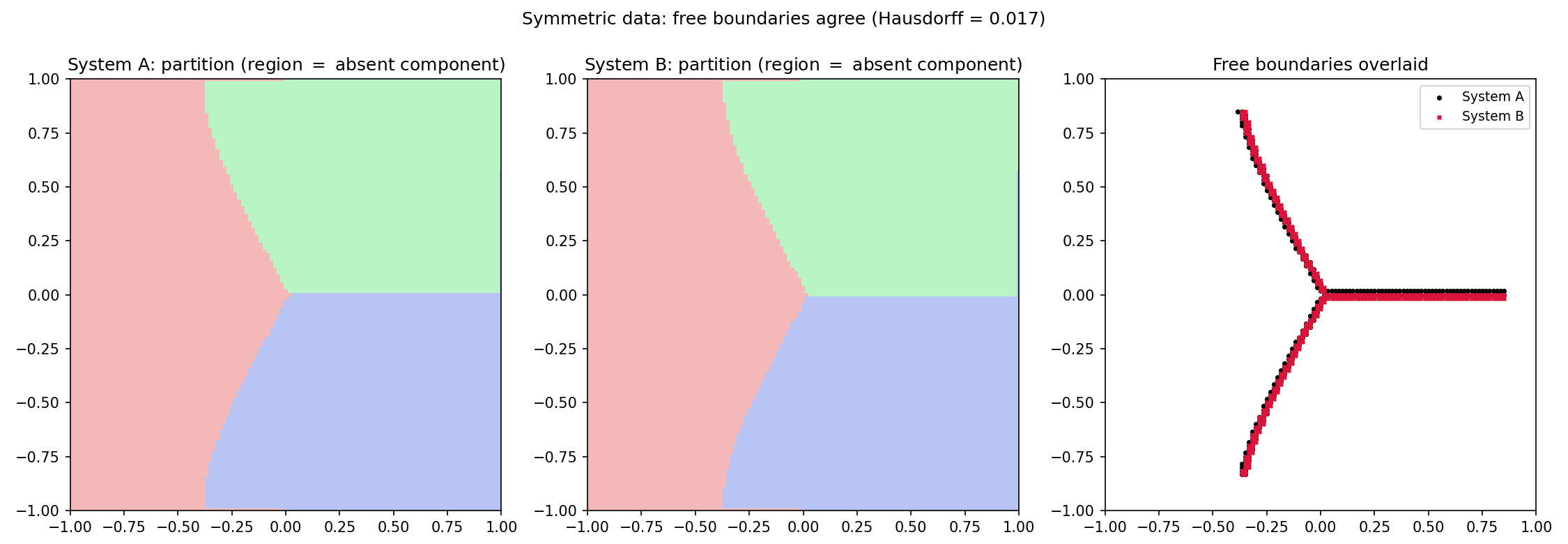}}{\fbox{\parbox{.88\textwidth}{\centering Symmetric numerical comparison.}}}
\caption{Datum~I. The discrete gradient triangle is nearly equilateral, and
the System~A and System~B proxy interfaces agree to approximately one mesh
width.}
\label{fig:symmetric-numerics}
\end{figure}

\begin{figure}[h]
\centering
\IfFileExists{I7.png}{\includegraphics[width=.9\textwidth]{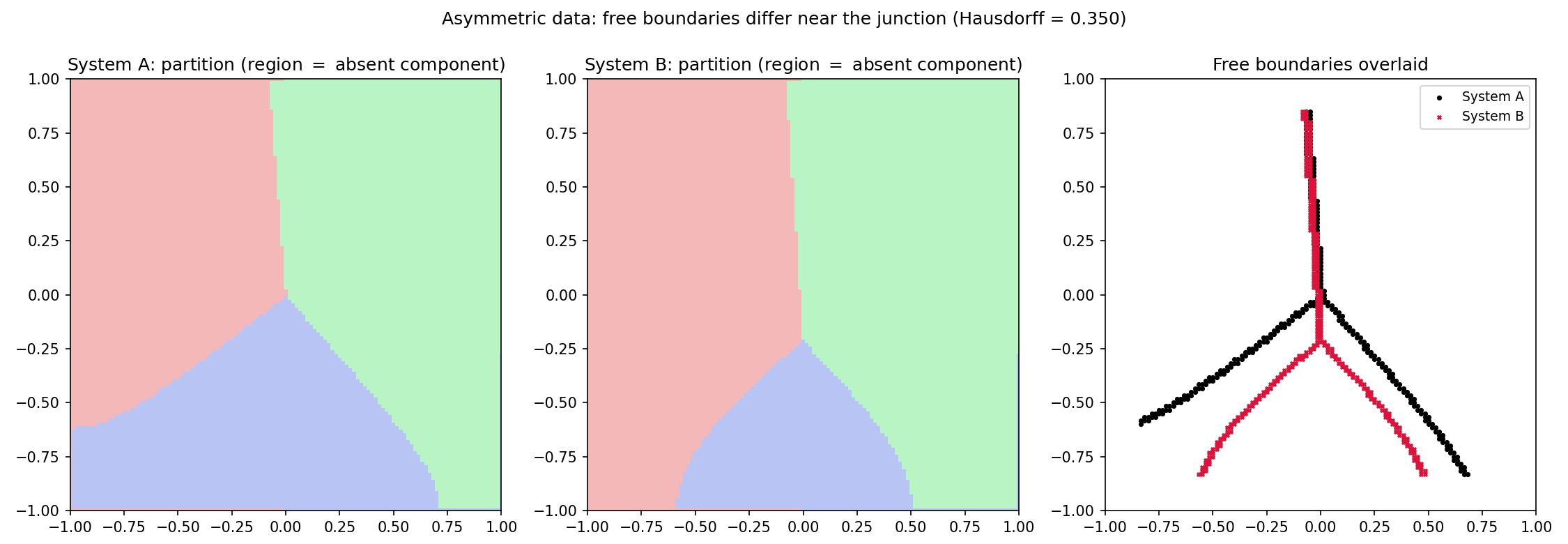}}{\fbox{\parbox{.88\textwidth}{\centering Asymmetric numerical comparison.}}}
\caption{Datum~II. The discrete gradient triangle is strongly
nonequilateral, and the finite-\(\eps\) System~B proxy is visibly displaced
from the System~A partition.}
\label{fig:asymmetric-numerics}
\end{figure}

The geometric diagnostics are summarized in Table~\ref{tab:gradtri}. Here
\(p_h\) is the approximate System~A triple tie, the quantities
\(|\nabla_h h_{ij}(p_h)|\) are the side lengths of the discrete gradient
triangle, and \(d_H\) is the symmetric Hausdorff distance between the two
proxy interfaces.

\begin{table}[h]
\centering
\caption{Geometric diagnostics on the \(121\times121\) grid.}
\label{tab:gradtri}
\renewcommand{\arraystretch}{1.25}
\begin{tabular}{lcccc}
\hline
Datum
& \(p_h\)
& \(\bigl(|\nabla_h h_{12}|,|\nabla_h h_{23}|,
|\nabla_h h_{31}|\bigr)(p_h)\)
& ray separations
& \(d_H\)\\
\hline
Datum~I
& \((0.018,\,0.000)\)
& \((0.807,\,0.791,\,0.807)\)
& \((121.4^\circ,\,119.3^\circ,\,119.3^\circ)\)
& \(0.017\)\\
Datum~II
& \((0.002,\,-0.007)\)
& \((0.557,\,0.773,\,0.922)\)
& \((123.2^\circ,\,93.8^\circ,\,142.9^\circ)\)
& \(0.358\)\\
\hline
\end{tabular}
\end{table}

\section{Conclusions and open problems}\label{sec:discussion}

The two limiting systems are governed by different selection mechanisms.
System~A is determined by global harmonic-difference identities and admits
the explicit lower-envelope representation, whereas System~B is selected by
constrained Dirichlet-energy minimization. Consequently, equality of their
interfaces is not automatic. The results above provide sufficient coincidence
criteria, an explicit one-dimensional noncoincidence example, and a necessary
gradient-triangle condition at common simple planar triple points.



\begin{thebibliography}{99}

\bibitem{BBF}
F. Bozorgnia, M. Burger, and M. Fotouhi, On a class of singularly perturbed elliptic systems with
asymptotic phase segregation, \emph{Discrete Contin. Dyn. Syst.} \textbf{42} (2022), 3539--3556.

\bibitem{BozorgniaArakelyanGamma}
F. Bozorgnia and A. Arakelyan,
Gamma convergence of partially segregated elliptic systems,
preprint, arXiv:2510.03794, 2025.

\bibitem{CaffarelliLin}
L. A. Caffarelli and F.-H. Lin, Singularly perturbed elliptic systems and multi-valued harmonic
functions with free boundaries, \emph{J. Amer. Math. Soc.} \textbf{21} (2008), 847--862.

\bibitem{CaffarelliRoquejoffre}
L. A. Caffarelli and J.-M. Roquejoffre, Uniform H\"older estimates in a class of elliptic systems and
applications to singular limits in models for diffusion flames, \emph{Arch. Ration. Mech. Anal.}
\textbf{183} (2007), 457--487.

\bibitem{ContiTerraciniVerzini}
M. Conti, S. Terracini, and G. Verzini, An optimal partition problem related to nonlinear
eigenvalues, \emph{J. Funct. Anal.} \textbf{198} (2003), 160--196.

\bibitem{GT}
D. Gilbarg and N. S. Trudinger, \emph{Elliptic Partial Differential Equations of Second Order},
Classics in Mathematics, Springer-Verlag, Berlin, 2001.

\bibitem{SoaveTerracini1}
N. Soave and S. Terracini, On some singularly perturbed elliptic systems modeling partial
segregation: uniform H\"older estimates and basic properties of the limits, preprint,
arXiv:2409.11976 (2024).

\bibitem{SoaveTerracini2}
N. Soave and S. Terracini, On partially segregated harmonic maps: optimal regularity and structure of
the free boundary, preprint, arXiv:2410.23976 (2024).

\bibitem{TavaresTerracini}
H. Tavares and S. Terracini, Regularity of the nodal set of segregated critical configurations under
a weak reflection law, \emph{Calc. Var. Partial Differential Equations} \textbf{45} (2012), 273--317.

\bibitem{Giaretto}
L. Giaretto,
On elliptic systems with \(k\)-wise interactions in the strong competition
regime: uniform H\"older bounds and properties of the limiting configurations,
preprint, arXiv:2603.10949, 2026.

\bibitem{OgnibeneVelichkovTriple}
R. Ognibene and B. Velichkov,
Structure of the free interfaces near triple junction singularities in
harmonic maps and optimal partition problems,
preprint, arXiv:2412.00781, 2024.

\bibitem{WangZhang}
K. Wang and Z. Zhang, Some new results in competing systems with many species, \emph{Ann. Inst. H.
Poincar\'e Anal. Non Lin\'eaire} \textbf{27} (2010), 739--761.

\end{thebibliography}
\end{document}